\documentclass[]{amsart}
\usepackage{pinlabel}

\usepackage{graphicx}
\usepackage{multirow}
\usepackage[english]{babel}
\usepackage{amsrefs}
\graphicspath{ {images/} }
\usepackage[T1]{fontenc}
\usepackage{txfonts}
\usepackage{times}
\usepackage[all]{xy}
\usepackage{subfig}
\usepackage{tikz}
\usetikzlibrary{shapes,arrows,shadows}
\usetikzlibrary{decorations.markings}
\usepackage{color}
\usepackage[normalem]{ulem}
\usepackage{hyperref}
\usepackage{wrapfig}
\usetikzlibrary{arrows}
\usepackage{pb-diagram}
\usepackage{verbatim}
\usepackage{float}
\usepackage{pstricks}
\usepackage{pst-plot}
\usepackage{graphpap}
\usepackage{pst-all}
\newtheorem{theorem}{Theorem}
\newtheorem{lemma}{Lemma}[section]
\newtheorem{proposition}{Proposition}
\newtheorem{corollary}{Corollary}
\newtheorem{definition}{Definition}

\newtheorem{remark}{Remark}
\newtheorem{claim}{Claim}
\newtheorem{construction}{Construction}[section]

\begin{document}

\title{On the Ends of groups and the Veech groups of infinite-genus surfaces}

\author[Ram\'irez Maluendas]{Camilo Ram\'irez Maluendas}
\address{Camilo Ram\'irez Maluedas\\\newline Departamento de Matem\'aticas y Estad\'istica, Universidad Nacional de Colombia, Sede Manizales, Manizales C.P. 170001, Colombia.}
\email{camramirezma@unal.edu.co}
\thanks{Partially supported by UNIVERSIDAD NACIONAL DE COLOMBIA, SEDE MANIZALES, Project Hermes 55598.}


\begin{abstract}
	
In this paper, we study the PSV construction, which provides a step by step method for obtaining tame translation surfaces with a suitable Veech group. In addition, we modify slightly this construction, and for each finitely generated subgroup $G<{\rm GL}_{+}(2,\mathbb{R})$ without contracting elements, we produce a tame translation surface $S$ with infinite genus such that its Veech group is $G$. Furthermore, the ends space of $S$ can be written as $\mathcal{B}\sqcup \mathcal{U}$, where $\mathcal{B}$ is homeomorphic to the ends space of the group $G$, and $\mathcal{U}$ is a countable, discrete, dense, and open subset of the ends space of $S$.
\end{abstract}

\maketitle

\section*{Introduction}\label{sec:introduction}

Geometrically, an \emph{end} of a topological space is a point at infinity. In \cite{Fre}, Freudenthal introduced the concept of ends and explored some of its applications in group theory. One can define the ends space ${\rm Ends}(G)$ of a finitely generated group $G$ as the ends space of the Cayley graph ${\rm Cay}(G,H)$, where $H$ is a generating set of $G$ (see \cite
{Fre1, Hopf}). In the context of orientable surfaces, Ker\'ekj\'art\'o \cite{Ker} studied their ends and introduced the classification of non-compact orientable surfaces, which determines the topological type of any orientable surface $S$ by its genus $g(S)\in \mathbb{N}\cup \{\infty\}$ and two closed subsets, ${\rm Ends}_{\infty}(S)\subseteq {\rm Ends}(S)$, of the Cantor set. These subsets are referred to as the ends space of $S$, and the ends of $S$ having (infinite) genus (see \cite{Ian}). Our focus is on studying surfaces with infinite genus.

\emph{Translation surfaces} have naturally appeared in various contexts: Dynamical systems (see \cite{KMS, HuSch}), Teichm\"uller theory (see \cite{KoZo, Moll}), Riemann surfaces (see \cite{MasTab, Zor}), among others. Our focus is on the so-called \emph{tame} translation surfaces. Using the charts of a translation surface $S$, one can pull back the standard Riemannian metric on $\mathbb{R}^2$ to equip the surface $S$ with a flat Riemannian metric $\mu$. This flat metric induces a distance map $d$ on $S$. A translation surface $S$ is said to be \emph{tame} \cite{Val2012} if, for each point $x\in \widehat{S}$ (where $\widehat{S}$ is the metric completion of $S$ with respect to $d$,) there is a neighborhood $U_{x}\subset \widehat{S}$ that is isometric to either an open subset of the Euclidean plane or, an open subset around a ramification point of a (finite or infinite) cyclic branched covering of the unit disk. It is worth noting that if $S$ is a compact translation surface, then $S$ is necessarily tame. Several authors have studied such surfaces (see for instance \cite{BowVal, FoKe, FrCo, Anja, VaWe}), which provides strong motivation for our research.

During the 1989s, Veech \cite{Vee} associated a group of matrices $\Gamma < {\rm GL}(2,\mathbb{R})$ to each translation surface, now commonly known as the \emph{Veech group of $S$}. He proved that if the Veech group $\Gamma(S)$ of a compact translation surface $S$ is a lattice--meaning $\Gamma(S)$ is a Fuchsian group such that the quotient space $\mathbb{H}^{2}/ \Gamma$ has finite hyperbolic area-- then the behavior of the geodesic flow on $S$ exhibits dynamical properties similar to those described by Weyl's theorem for the geodesic flow on the torus. This result is known as the \emph{Veech's Dichotomy}. It has since attracted the attention of many researchers (see, for example, \cite{Fin, Hoop, PasLan}).

The Veech group associated to a compact translation surface is a Fuchsian group \cite{Voro}. In the case of a tame translation surface,  if $\Gamma(S)$ is the Veech group of the tame translation surface $S$, then one of the following holds \cite[Theorem 1.1]{PSV}:
\begin{enumerate}
	\item[\textbf{(1)}] $\Gamma(S)$ is countable and without contracting elements, it means,  $\Gamma(S)$ is disjoint from the set $\{A\in {\rm GL}_{+}(2,\mathbb{R}): \Vert Av\Vert < \Vert v\Vert  \text{ for all } v\in\mathbb{R}^2\setminus \{\textbf{0}\}\}$, where $\Vert \, \Vert$ is the Euclidean norm on $\mathbb{R}^2$, or
	\item[\textbf{(2)}] $\Gamma(S)$ is conjugated to
	$
	P:=\left\{
	\begin{pmatrix}
		1 & t \\
		0 & s
	\end{pmatrix}\hspace{1mm}:\hspace{1mm} t\in\mathbb{R} \text{  and }
	\hspace{1mm} s\in\mathbb{R}^{+}
	\right\},\hspace{1mm}\text{or}
	$
	\item[\textbf{(3)}] $\Gamma(S)$ is conjugated to $P'<{\rm GL_{+}(2, \mathbb{R})}$, the subgroup generated by $P$ and $\rm -Id$, or
	\item[\textbf{(4)}] $\Gamma(S)$ is equal to ${\rm GL_{+}(2, \mathbb{R})}$.
\end{enumerate}

Our work contributes to the problem of realizing subgroups of ${\rm GL}_{+}(2,\mathbb{R})$ as Veech groups of (non-compact) tame translation surfaces. We will discuss some of the studies involved in the problem of realizing groups as symmetry groups of translation surface. In \cite{PSV}, the authors developed a step-by-step process, referred to as the \emph{PSV construction}, aimed at constructing, for each subgroup $G < {\rm GL}(2,\mathbb{R})$ without contracting elements, a tame Loch Ness monster with Veech group $G$. Up to homeomorphism, the \emph{Loch Ness monster} is the only surface with infinite genus and a unique end \cite{PSul}. In the case of \emph{origamis}, translation surfaces formed by appropriately gluing unit squares, any finite group can be represented as the automorphism group of the Loch Ness monster when it is viewed as an origami \cite{HIMO}. The PSV construction, with slight modifications, was used in \cite{RamVal} to realize any subgroup $G < {\rm GL}_{+}(2,\mathbb{R})$ without contracting elements as the Veech group of a large class of tame translation surfaces of infinite genus. These results, along with those addressing the realization of Veech groups for translation surfaces with non-self-similar end spaces \cite{MoraVal}, have been extended to resolve the problem of realizing symmetry groups of infinite genus translation surfaces \cite{ARSVW}.

We have also explored and made slightly modifications to the PSV construction, resulting in a theorem that establishes an explicit connection between the ends space of a tame translation surface and the ends space of its respective Veech group.

%%%%%%%%%%%%%%%%%%%%%%%%%%%%%%%%%%%%%%%%%%%%%%%%%%%%%%%%%

\begin{theorem}\label{theorem.0.2}
	Given a finitely generated subgroup $G$ of ${\rm GL}_{+}(2,\mathbb{R})$ without contracting elements. Then there exists a tame translation surface $S$ whose Veech group is $G$. The ends space ${\rm Ends}(S)$ of $S$ satisfies:
	
	\begin{enumerate}
		\item[\textbf{(1)}] If $G$ is finite, then the surface $S$ has as many ends as there are elements in the group $G$, and each end has infinite genus.
		
		\item[\textbf{(2)}]  If $G$ is not finite, then the ends space of $S$  can be represented as  
		\[
		{\rm Ends}(S)={\rm Ends}_{\infty}(S)=\mathcal{B}\sqcup \mathcal{U},
		\]
		where $\mathcal{B}$ is a closed subset of ${\rm Ends}(S)$ homeomorphic to ${\rm Ends}(G)$, and $\mathcal{U}$ is a countable, discrete, dense, and open subset of ${\rm Ends}(S)$.
	\end{enumerate}
\end{theorem}

As the ends space of a finitely generated group has either zero, one, two, or infinitely many ends \cite{Fre1, Hopf}, we immediately obtain the following corollary: 

\begin{corollary}
	The ends space of the tame translation surface $S$ is one of the following:
	
	\begin{enumerate}
		\item[\textbf{(1)}] If the group $G$ has one end, then ${\rm Ends}(S)$ is homeomorphic to the ordinal number $\omega+1$. In other words, the ends space of $S$ is homeomorphic to the closure of $\left\{\frac{1}{n}:n\in\mathbb{N}\right\}$.
		
		\item[\textbf{(2)}] If the group $G$ has two ends, then ${\rm Ends}(S)$ is homeomorphic to the ordinal number $\omega\cdot 2+1$. This means that the ends space of $S$ is homeomorphic to two copies of the closure of $\left\{\frac{1}{n}:n\in\mathbb{N}\right\}$.

		\item[\textbf{(3)}] If the group $G$ has infinitely many ends, then ${\rm Ends}(S)$ contains a subset homeomorphic to the Cantor set, with its complement being a countable, discrete, dense, and open subset of ${\rm Ends}(S)$.
	\end{enumerate}
\end{corollary}

The paper is structured as follows: In Section \ref{sec:ends}, we collect the principal tools needed to understand the classification of non-compact surfaces theorem and explore the concept of ends on groups. Section \ref{section:tame_veech} provides an introduction to the theory of tame translation surfaces and discusses the Veech group. Finally, Section \ref{proof} is dedicated to proving our main result.

%%%%%%%%%%%%%%%%%%%%%%%%%%%%%%%%%%%%%%%%%%%%%%%%%%
%%%%%%%%%%%%%%%%%%%%%%%%%%%%%%%%%%%%%%%%%%%%%%%%%%

%%%%%%%%%%%%%%%%%%%%%%%%%%%%%%%%%%%%%%%%%%%%%%%%%
%%%%%%%%%%%%%%%%%%%%%%%%%%%%%%%%%%%%%%%%%%%%%%%%%

\section{Ends}\label{sec:ends}

In this section, we shall introduce the concept of the space of ends of a topological space $X$ in its most general context. We shall also explore the classification theorem of non-compact orientable surfaces based on their ends spaces.  Finally, we shall discuss the concept of ends on groups. 

\begin{definition}[\cite{Fre}] \label{section_ends_space} 
	Let $X$ be a locally compact, locally connected, connected, and Hausdorff space, and let $(U_n)_{n\in\mathbb{N}}$ be  an infinite nested sequence $U_{1}\supset U_{2}\supset\ldots$ of non-empty connected open subsets of $X$, such that the following conditions hold:
	\begin{enumerate}
		\item[\textbf{(1)}]  For each $n\in\mathbb{N}$, the boundary $\partial U_n$ of $U_{n}$ is compact. 
		
		\item[\textbf{(2)}] The intersection $\bigcap\limits_{n\in\mathbb{N}}\overline{U_{n}}=\emptyset
		$.
		
		\item[\textbf{(3)}] For any compact  subset $K\subset X$, there is $m\in\mathbb{N}$ such that $K\cap U_m =\emptyset$.
	\end{enumerate}
	Two nested sequences $(U_n)_{n\in\mathbb{N}}$ and $(U'_{n})_{n\in\mathbb{N}}$ are {\it equivalent} if, for each $n\in\mathbb{N}$, there exist $j,k\in\mathbb{N}$ such that $U_{n}\supset U'_{j}$, and $U'_{n}\supset U_{k}$. The corresponding equivalence classes of these sequences are called the \textbf{ends} of $X$. The ends space ${\rm Ends}(X)$ of $X$ is the space whose elements are the ends of $X$, and it is endowed with the following topology: for any non-empty open subset $U$ of
	$X$, such that its boundary $\partial U$ is compact, we define
	\begin{equation}\label{eq:end_open}
		U^{*}:=\left\{[U_{n}]_{n\in\mathbb{N}}\in{\rm Ends}(X)\hspace{1mm}|\hspace{1mm}U_{j}\subset U\hspace{1mm}\text{for some }j\in\mathbb{N}\right\}.
	\end{equation}
	Then the set of all such $U^{\ast}$, where $U$ is open and has a compact boundary in $X$, forms a basis for the topology of ${\rm Ends}(X)$ (see \cite[1. Kapitel]{Fre}).
\end{definition}

\begin{theorem}[\cite{Ray}]
	The space ${\rm Ends}(X)$, with the topology defined above, is Hausdorff, totally disconnected,
	and compact.
\end{theorem}

\subsection{Ends of a surface}

A \emph{surface} S is a connected $2$-manifold without boundary, which may or may not be closed. In this manuscript, we shall only consider orientable surfaces. By a \emph{subsurface} of $S$ we mean an embedded surface, which is a closed subset of $S$, and whose boundary consists of a finite number of nonintersecting simple closed curves. Note that a subsurface may or may not be compact.
The {\it reduced genus} of a compact subsurface $\tilde{S}\subset S$,  with $q(\tilde{S})$ boundary curves and Euler characteristic $\chi(\tilde{S})$, is the number 
\[
g(\tilde{S})=1-\frac{1}{2}\left(\chi(\tilde{S})+q(\tilde{S})\right).
\] 
The {\it genus} of the surface $S$ is the supremum of the genera of its compact subsurfaces. This genus may be a non-negative integer or $\infty$. 
The surface $S$ is said to be {\it planar} if it has genus zero, in other words, $S$ is homeomorphic to an open of the complex plane. 

\begin{remark}
	In this case, from the definition of ends given in Definition \ref{section_ends_space}, we may assume that for the sequence $(U_{n})_{n\in\mathbb{N}}$ the closures $\overline{U}_{n}$ are subsurfaces. In this setting, an end $[U_n]_{n\in\mathbb{N}}$ of a surface $S$ is called \textbf{planar}, if there is $l\in\mathbb{N}$ such that the subsurface $\overline{U}_l\subset S$ is planar.
\end{remark}

We define the subset ${\rm Ends}_{\infty}(S)$ of ${\rm Ends}(S)$ to co consist of all ends of $S$, which are not planar ({\it ends having infinity genus}). It follows directly from the definition that ${\rm Ends}_{\infty}(S)$ is a closed subset of ${\rm Ends}(S)$ (see \cite[p. 261]{Ian}), and the triplet $(g,{\rm Ends}_{\infty}(S),{\rm Ends}(S))$, where $g$ is the genus of $S$, is a topological invariant.

\begin{theorem}[Classification of non-compact surfaces \cite{Ker, Ian}]
	Two surfaces $S_1$ and $S_2$  having the same genus, are topologically equivalent if and only if there exists a homeomorphism $f: {\rm Ends}(S_1)\to {\rm Ends}(S_2)$ such that $f( {\rm Ends}_{\infty}(S_1))= {\rm Ends}_{\infty}(S_2)$.
\end{theorem}

\begin{definition}[\cite{PSul}]
	The \textbf{Loch Ness monster} is the unique, up to homeomorphism, infinite genus surface with exactly one end. 
\end{definition}

\begin{remark}[\cite{SPE}] The surface $S$ has $m$ ends, for some $m\in\mathbb{N}$, if and only if for any compact subset $K \subset S$, there is a compact $K^{'}\subset S$ such that $K\subset K^{'}$ and $S\setminus  K^{'}$ consists of $m$ connected components.
\end{remark}

\subsection{Ends of a group}

Given a generating set $H$ (closed under inverse) of a group $G$, the {\it Cayley graph of $G$ with respect to the generating set $H$} is the graph ${\rm Cay}(G, H)$, where the vertices are the elements of $G$, and there is an edge between two vertices $g_1$ and $g_2$ if and only if there is $h\in H$ such that $g_1h=g_2$. Throughout this paper, the Cayley graph ${\rm Cay}(G, H)$ will be the geometric realization of an abstract graph \cite[p. 226]{Diestel2}.

When the set $H$ is finite, the Cayley graph ${\rm Cay}(G, H)$ is locally compact, locally connected, connected, and Hausdorff space. In this case, we define the {\it ends space of $G$} as ${\rm Ends}(G):={\rm Ends}({\rm Cay}(G, H))$.

\begin{proposition}[\cite{Loh}] 
	Let $G$ be a finitely generated group. The ends space of the Cayley graph of $G$ does not depend on the choice of the finite generating set.
\end{proposition}

\begin{theorem}[\cite{Fre1, Hopf}]
	Let $G$ be a finitely generated group. Then $G$ has either zero, one, two, or infinitely many ends.
\end{theorem}

%%%%%%%%%%%%%%%%%%%%%%%%%%%%%%%%%%%%%%%%%%%%%%%%%%%%%%%%
%%%%%%%%%%%%%%%%%%%%%%%%%%%%%%%%%%%%%%%%%%%%%%%%%%%%%%%%

\section{Tame translation surfaces}\label{section:tame_veech}

An atlas $\mathcal{A}=\{(U_{\alpha},\phi_{\alpha})\}_{\alpha\in I}$ on the  surface $S$ is called a \emph{translation atlas} if $S$, except for a subset of points ${\rm Sing}(S)\subset S$, can be covered by the charts from such atlas. Moreover, for any pair of charts $(U_{\alpha},\phi_{\alpha})$ and $(U_{\beta},\phi_{\beta})$ in $\mathcal{A}$ such that $U_{\alpha} \cap U_{\beta}\neq \emptyset$, the associated transition map
\[
\phi_{\alpha}\circ \phi_{\beta}^{-1}:\phi_{\beta}(U_{\alpha} \cap U_{\beta})\subset \mathbb{R}^2\to \phi_{\alpha}(U_{\alpha} \cap U_{\beta})\subset\mathbb{R}^2,
\]
is locally the restriction of a translation. We assume that each point in ${\rm Sing}(S)$ is non-removable, which means, the translation atlas can not be extended to any of the points in ${\rm Sing}(S)$. An element $x$ in ${\rm Sing}(S)$ is called a \emph{singular point of $S$} or \emph{singularity}. A \emph{translation structure} on $S$ is a maximal translation atlas on the surface. If $S$ admits a translation structure, it will be called a \emph{translation surface}. 

For a translation surface $S$, we can pull back the Euclidean (Riemannian) metric of $\mathbb{R}^{2}$ via its translation structure, thus we obtain a flat Riemannian metric $\mu$ on $S$. Let $\widehat{S}$ denote the \emph{metric completion of $S$} with respect to the flat Riemannian metric $\mu$. According to the Uniformization Theorem \cite[p. 580]{Abi}, the only complete translation surfaces  $S=\widehat{S}$ are the Euclidean plane, the torus, and the cylinder \cite[p. 193]{FarKra}.

\begin{definition}[\cite{Val2012}]
	A translation surface $S$ is said to be  \textbf{tame} if for each point $x\in \widehat{S}$, there exists a neighborhood $U_x \subset \widehat{S}$ isometric to either:
	\begin{itemize}
		\item[\textbf{(1)}] Some open subset of the Euclidean plane, or
		
		\item[\textbf{(2)}] An open subset of the ramification point of a (finite or infinite) cyclic branched covering of the unit disk in the Euclidean plane.
	\end{itemize}
	
	In the second case, if the neighborhood $U_{x}$ is isometric to the finite cyclic branched covering of finite order $m\in\mathbb{N}$, then the point $x$ is called a \textbf{finite cone angle singularity of angle $2m\pi$}. If $U_{x}$ is isometric to the infinite cyclic branched covering, then $x$ is called a \textbf{infinite cone angle singularity}. 
\end{definition}

We denote by ${\rm Sing}(\widehat{S})$ the set of all the finite and infinite cone angle singularities of $\widehat{S}$. An element of ${\rm Sing}(\widehat{S})$ is called a \emph{cone angle singularity of} $\widehat{S}$, or simply a \emph{cone point}.

\subsection{Saddle connection and markings}

A \emph{saddle connection} $\gamma$ on a tame translation surface $S$ is a geodesic interval joining two cone points and not having cone points in its interior. In the translation structure of $S$, we can find a chart $(U,\varphi)$ such that the open $U$ contains the saddle connection $\gamma$, excluding its endpoints. The map $\varphi$ sends $\gamma$ to a straight line segment in $\mathbb{R}^{2}$. This straight line segment can be oriented in two possible directions, denoted by $[\theta],[-\theta]\in \mathbb{R}/2\pi \mathbb{Z}$, for some $\theta \in\mathbb{R}$. Then we can associate to $\gamma$ two oppositely oriented vectors $\{v,-v\}\subset\mathbb{R}^2$, corresponding to the direction $[\theta]$ and $[-\theta]$, respectively. Moreover, the norm of these vectors is equal to the length of $\gamma$, measured with respect to the flat Riemannian metric $\mu$ on $S$. Each one of these vectors is called a \emph{holonomy vector of $\gamma$}. Clearly, the holonomy vectors of $\gamma$ are well-defined, that is, does not depend on the choice of the chart $(U,\varphi)$.

A \emph{marking} $m$ on the tame translation surface $S$ is a finite length geodesic not having cone points inside it. Similarly to the case of saddle connection, we can associate to the marking $m$ two \emph{holonomy vectors} $\{v,-v\}\subset \mathbb{R}^2$. Two markings are said to be \emph{parallel} if their respective holonomy vectors are also parallel. It does not matter if the markings are on different surfaces \cite[Definition 3.4]{PSV}. 

\begin{definition}[\cite{RamVal}]
	Let $m_1$ and $m_2$ be two parallel markings having the same length on translation surfaces $S_1$ and $S_2$, respectively. We cut $S_1$ and $S_2$ along $m_1$ and $m_2$, respectively, turning $S_1$ and $S_2$ into the surfaces with boundary $\tilde{S}_1$ and $\tilde{S}_2$, respectively. Each one of their boundaries is formed by two straight line segments. Now, we consider the union $\tilde{S}_1\cup \tilde{S}_2$ and identify  (glue) such (four) segments using translations to obtain a connected tame translation surface $S$ (see Figure \ref{Figure_gluing_marks}). This gluing relation of these segments will be denoted as $m_1\sim_{\text{glue}}m_2$, and will be called \textbf{the operation of gluing the markings $m_1$ and $m_2$}. Then the surface $S$ will be written in the form 
	\[
	S:=(S_1\cup S_2)/m_1\sim_{\text{glue}}m_2.
	\]
	We say that $S$ is obtained from $S_1$ and $S_2$ by \textbf{regluing} along $m_1$ and $m_2$.

	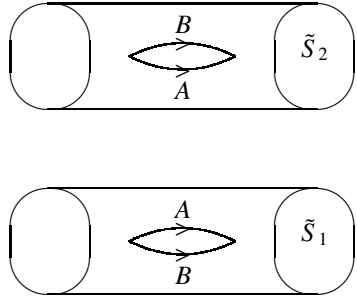
\begin{figure}[ht!]
		\begin{center}
			\begin{tikzpicture}[baseline=(current bounding box.north), scale=0.8]
				\begin{scope}
					\clip (-4,-0.4) rectangle (13,5.2);
					%%%%%%%%%%%%%%%%%%%%%%%%%%%%%%%%%%%%%%%%%%%%%%%
					%%
					\draw [line width=1pt]  (2,1) ellipse (5mm and 10mm);
					\draw [line width=1pt]  (2,4) ellipse (5mm and 10mm);
					\draw [line width=1pt]  (6,1) ellipse (5mm and 10mm);
					\draw [line width=1pt]  (6,4) ellipse (5mm and 10mm);
					\draw [line width=1pt]  (2,0) -- (6,0);
					\draw [line width=1pt]  (2,2) -- (6,2);
					\draw [line width=1pt]  (2,3) -- (6,3);
					\draw [line width=1pt]  (2,5) -- (6,5);
					%%%%%%%%%%%%%%%%%%%%%%%%%%%%%%%%%%%%%%
					\draw [->, color=red, line width=1pt] plot[smooth] coordinates
					{(3,1)(4,0.7)(5,1)};
					\draw [->, color=blue, line width=1pt] plot[smooth] coordinates
					{(5,1)(4,1.3)(3,1)};
					%%%%%%%%%%%%%%%%%%%%%%%%%%%%%%%%%%%%%%%%%%%%
					\draw [<-, color=blue, line width=1pt] plot[smooth] coordinates
					{(3,4)(4,3.7)(5,4)};
					\draw [<-, color=red, line width=1pt] plot[smooth] coordinates
					{(5,4)(4,4.3)(3,4)};
					%%%%%%%%%%%%%%%%%%%%%%%%%%%%%%%%
					%%%%%%%%%%%%%%%%%%%%%%%%%%%%%%%%%%%%%%%
					%%%%%%%%%%%%%%%%%%%%%%%%%%%%%%
					%%%%%%%%%%%%%%%%%%%%%%%%%%%%%%%%%%%%%%%%%
					\node at (4,0.4) {$B$};
					\node at (4,1.6) {$A$};
					\node at (4,3.4) {$A$};
					\node at (4,4.6) {$B$};
				\end{scope}
			\end{tikzpicture}
		\end{center}
		\caption{\emph{Gluing markings.}}
		\label{Figure_gluing_marks}
	\end{figure}
\end{definition}

\subsection{Veech group}
Let $S$ be a tame translation surface. A homeomorphism $T:\widehat{S}\to \widehat{S}$ is called \emph{affine diffeomorphism}, if it satisfies the following properties:
\begin{itemize}
	\item[\textbf{(1)}] It sends cone points to cone points.
	\item[\textbf{(2)}] The function $T$ is an affine map in the local coordinates of the translation atlas on $S$.
\end{itemize}
We denote by ${\rm Aff}_{+}(S)$ the group of all the affine orientation preserving diffeomorphism from the tame translation surface $S$ to itself. 

Given a tame translation surface $S$ and a map $T\in {\rm Aff}_{+}(S)$, then using the translation structure on $S$, we hold that the differential $dT(p)$ of $T$ at any point $p\in S$ is a constant matrix $A$ belongs to ${\rm GL}_{+}(2,\mathbb{R})$. We then define the map
\begin{equation}\label{eq:homomorphism_differential_matrix}
	D: {\rm Aff}_{+}(S)\to {\rm GL}_{+}(2,\mathbb{R}),
\end{equation}
where $D(T)$ is the differential matrix of $T$. Using the chain rule, it is easy to verify that $D$ is a group homomorphism. 

\begin{definition}[\cite{Vee}]\label{def:Veech_group} 
	The image of $D$, denoted by $\Gamma(S)$, is called the \textbf{Veech group of $S$}. 
\end{definition}

The group ${\rm GL}_{+}(2,\mathbb{R})$ acts on the set of all translation surfaces by post-composition on charts. More precisely, this action sends the couple $(g, S)$ to the translation surface $S_g$, which is called \emph{the affine copy of $S$}. The translation structure on $S_{g}$ is obtained by postcomposing each chart on $S$ by the affine transformation associated to the matrix $g$. Further, this action defines an affine diffeomorphism $f_g:S\to S_{g}$, where the differential $d f_{g}(p)$ of $f_{g}$ at any point $p\in S$ is the matrix $g$.

\section{Proof of Theorem \ref{theorem.0.2}}\label{proof}

Let $G$ be a finitely generated subgroup of ${\rm GL}_{+}(2,\mathbb{R})$ without contracting elements, and let $H$ be a finite generating set of $G$. The set $H$ can be written as $H=\{h_j: j\in\{1,\ldots, J\}\}$, for some $J\in\mathbb{N}$. We shall obtain the surface $S$ using the PSV construction, which will be briefly outlined below. Afterward, we shall prove that $S$ is a tame translation surface with Veech group $G$. Finally, we will describe the ends space of $S$.

\subsection{PSV construction}

For each countable subgroup $G$ of ${\rm GL}_{+}(2,\mathbb{R})$ without contracting elements, Przytycki, Weitze-Schmith{\"u}sen, and Valdez, in \cite[4. Countable Veech group]{PSV}, described a method to construct a tame translation surface homeomorphic to the Loch Ness monster, with Veech group $G$. We refer to this method as the \emph{PSV construction}. From a metric spaces point of view, the process is as follows:

\subsubsection*{Step 1. The decorated surface}

We build a \emph{suitable} tame Loch Ness monster $S_{{\rm dec}}$ using copies of the Euclidean plane and a cyclic branched covering of the Euclidean plane, which are appropriately attached via gluing markings. The resulting surface $S_{\rm dec}$ is referred to as \emph{decorated}. For each $h_{j}\in H$, we mark $S_{{\rm dec}}$ with two infinite families of (suitable) markings 
\[
h_{j}\check{M}^{-j}:=\left\{h_{j}\check{m}^{-j}_{i}: \forall i\in\mathbb{N}\right\} \text{ and } M^{-j}:=\left\{m^{-j}_{i}: \forall i\in\mathbb{N}\right\}.
\]

\subsubsection*{Step 2. The puzzle associated to de triplet $(1,G,H)$}

For each $g\in G$, we take the affine copy $S_g$ of the decorated surface $S_{{\rm dec}}$. We then define two families of markings on $S_{g}$: 
\[
gh_{j}\check{M}^{-j}:=\left\{gh_{j}\check{m}^{-j}_{i}: \forall i\in\mathbb{N}\right\} \text{ and } gM^{-j}:=\left\{gm^{-j}_{i}: \forall i\in\mathbb{N}\right\}, 
\]
These families corresponded to the image of $h_{j}\check{M}^{-j}$ and $M^{-j}$ on $S_{\rm dec}$ (respectively) under the diffeomorphism $f_{g}: S_{\rm dec}\to S_{g}$. Thus, we define the \emph{puzzle associated to the triplet $(1,G,H)$} as 
\[
\mathfrak{P}(1,G,H):=\left\{S_g: g\in G\right\},
\]
as is defined in \cite[Definition 3.1]{RamVal}. The term $1$ means that the decorated surface has only one end.

\textbf{Step 3. The assembled surface $S$ to the puzzle $\mathfrak{P}(1,G,H)$.}

We define the \emph{assembled surface to the puzzle $\mathfrak{P}(1,G,H)$} (see \cite[Definition 3.1]{RamVal}) as follows:
\begin{equation*}\label{assembled_surface_1}
	S:= \bigcup\limits_{g\in G}S_g \bigg/ \sim,
\end{equation*}
where $\sim$ is the equivalent relation given by the following gluing of the markings: for each edge $(g,gh_{j})$ of the Cayley graph ${\rm Cay}(G, H)$, the marking
$gh_{j}\check{m}^{-j}_{i}$ on $S_{g}$ is glued to the marking $gh_{j}m^{-j}_{i}$ on $S_{gh_{j}}$, for each $i\in\mathbb{N}$.

\subsection{We employ PSV construction to obtain the surface $S$}

\subsubsection*{Step 1. The decorated surface} 

The following auxiliary construction is necessary to obtain the decorated surface.  

\begin{construction}[Buffer surface]
	For each $j\in\{1,\ldots, J\}$, we consider $\mathbb{E}(j,1)$ and $\mathbb{E}(j,2)$ copies of the Euclidean plane, which are endowed with a fixed origin $\textbf{0}$ and an orthogonal basis $\beta=\{e_1,e_2\}$. We define markings on these surfaces, which are described by their endpoints. On $\mathbb{E}(j,1)$, we draw the families of markings: 
	\begin{align*}
		\check{M}^{j} &:= \left\{ \check{m}^{j}_{i} := (4ie_{1}, (4i+1)e_{1}) :\forall i \in \mathbb{N}\right\}, \text{ and }\\
		L &:= \left\{ l_{i} := ((4i+2)e_{1}, \, (4i+3)e_{1}):\forall i \in \mathbb{N}\right\}.
	\end{align*}
	
	On $\mathbb{E}(j,2)$ we take the family of markings: 
	\[
	L^{'} :=  \left\{l^{'}_{i} := ((2i+1)e_{2}, \, e_{1}+(2i+1)e_{2}) :\forall i \in \mathbb{N}\right\},
	\]
	and the marking: 
	\[
	h_{j}\check{m}^{-j}  :=  (2e_{2}, \, e_{1}+2e_{2}).
	\] 
	Finally, the marking  $l_i\in L$ on $\mathbb{E}(j,1)$ and the marking $l^{'}_{i}\in L^{'}$ on $ \mathbb{E}(j,2)$ are glued, for each $i\in\mathbb{N}$. Thus, we obtain a tame Loch Ness monster 
	\begin{equation}\label{buffer_surface}
		S(Id,h_j),
	\end{equation}
	which is called the \textbf{buffer surface associated to the element $h_j$ of $H$} (see Figure \ref{Figure_buffer}).
	\begin{figure}[ht!]
		\begin{center}
			\setlength{\unitlength}{0.8pt}
			\begin{picture}(350,180)
				%\graphpaper(0,0)(350,200)
				\put(90,0){\framebox(250,60)}
				\put(100,30){\vector(1,0){210}}
				\put(120,10){\vector(0,1){40}}
				\put(111,18){{\small$\textbf{0}$}}
				%%%%%%%%%%%%%%
				\put(140,30){\linethickness{0.7mm}\line(1,0){20}}
				\put(145,36){{\small$\check{m}_{1}^{j}$}}
				\put(180,30){\linethickness{0.7mm}\line(1,0){20}}
				\put(185,35){{\small $l_{1}$}}
				\put(220,30){\linethickness{0.7mm}\line(1,0){20}}
				\put(225,36){{\small$\check{m}_{2}^{j}$}}
				\put(260,30){\linethickness{0.7mm}\line(1,0){20}}
				\put(265,35){{\small$l_{2}$}}
				\put(290,40){$\ldots$}
				\put(290,10){{\small$\mathbb{E}(j,1)$}}
				%%%%%%%%%%%%%%%%
				%%%%%%%%%
				\put(0,0){\framebox(70,180)}
				\put(10,20){\vector(1,0){50}}
				\put(30,10){\vector(0,1){155}}
				\put(22,10){{\small$\textbf{0}$}}
				%%%%%%
				\put(30,40){\linethickness{0.7mm}\line(1,0){20}}
				\put(40,46){{\small$l_{1}^{'}$}}
				\put(30,60){\linethickness{0.7mm}\line(1,0){20}}
				\put(35,65){{\small$h_j\check{m}^{-1}$}}
				\put(30,85){\linethickness{0.7mm}\line(1,0){20}}
				\put(40,91){{\small$l_{2}^{'}$}}
				\put(30,130){\linethickness{0.7mm}\line(1,0){20}}
				\put(40,136){{\small$l_{3}^{'}$}}
				\put(40,155){$\vdots$}
				\put(2,168){{\small$\mathbb{E}(j,2)$}}
			\end{picture}
			\caption{\emph{Buffer surface $S(Id,h_j)$.}}
			\label{Figure_buffer}
		\end{center}
	\end{figure}
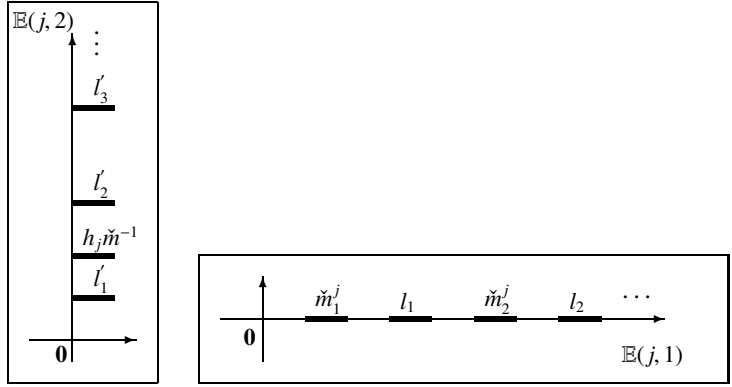 
\end{construction}

\begin{remark}
	The buffer surface $S(Id,h_j)$ is a modification of the surface appearing in Construction 4.4 in \cite{PSV}. We emphasize that the family of markings $\check{M}^{j}$ and the marking $h_{j}\check{m}^{-j}$ on $S(Id,h_j)$ have not yet been glued. In addition, the set of singular points of $S(Id, h_j)$ consists of infinitely many cone angle singularities of angle $4\pi$.
\end{remark} 

\begin{construction}[Decorated surface]
	We take $\mathbb{E}$, the Euclidean plane, endowed with a fixed origin $\overline{\textbf{0}}$ and an orthogonal basis $\beta= \{e_{1},e_{2}\}$. Analogously, we shall define markings on this surface, described by their endpoints. For each $j\in\{1,\ldots, J\}$, on $\mathbb{E}$ we define the families of markings: 
	\begin{align*}
		M^{j} &:=\left\{m^{j}_{i} := \left((2i-1)e_{1}+je_{2}, \, 2ie_{1}+je_{2}\right): \forall  i \in \mathbb{N} \right\}, \text{ and }\\ 
		M &:=\left\{m_i:=((4i-1)e_1,4ie_1): \forall i\in\mathbb{N}\right\}.
	\end{align*}
	
	Now, we recursively draw new markings on $\mathbb{E}$. For $j=1$, we choose two suitable real numbers $x_1>0$ and $y_1<0$, and we define the marking: 
	\[
	m^{-1} := (x_{1}e_{1}+y_{1}e_{2}, \, x_{1}e_1+h^{-1}_{1}e_{1}+y_{1}e_{2}),
	\]
	on $\mathbb{E}$, such that $m^{-1}$ is disjoint from the families of markings $M$ and $M^{j}$ for each $j\in\{1,\ldots, J\}$. 
	
	For $n\leq J$, we choose two suitable real numbers $x_n>0$ and $y_n<0$, and we define the marking: 
	\[
	m^{-n}:=(x_{n}e_{1}+y_{n}e_{2}, \, x_{n}e_1+h^{-1}_{n}e_{1}+y_{n}e_{2}),
	\]
	on $\mathbb{E}$, such that $m^{-n}$ is disjoint from the families of markings $M$ and $M^j$, for each $j\in \{1,\ldots, J\}$. Moreover, the marking $m^{-n}$ is also disjoint from the markings $m^{-1},\ldots, m^{-(n-1)}$ defined in the previous steps.

	Let $\pi: \tilde{\mathbb{E}}\to \mathbb{E}$ be the three fold cyclic covering of $\mathbb{E}$, branched over the origin. Then we denote as 
	\[
	\tilde{M}:=\{\tilde{m_{i}}: \forall i \in \mathbb{N}\},
	\]
	one of the three sets of markings on $\tilde{\mathbb{E}}$ defined by $\pi^{-1}(M)$. Now, we take on $\mathbb{E}$ the markings $t_{1}:= (e_{2}, \, 2e_{2})$, and $t_{2} := (-e_{2},\, -2e_{2})$, which will be used to generate new markings on $\tilde{\mathbb{E}}$. Then we denote as $\tilde{t_{1}}$ and $\tilde{t_{2}}$ one of the three markings on $\tilde{\mathbb{E}}$ defined by $\pi^{-1}(t_1)$ and $\pi^{-1}(t_2)$, respectively, such that they are on the same fold of $\tilde{\mathbb{E}}$ as $\tilde{M}$.
	
	Finally, we take the union of surfaces $\mathbb{E}\cup \tilde{\mathbb{E}}\bigcup_{j\in \{1,\ldots, J\}}S(Id,h_{j})$ (see equation \eqref{buffer_surface}), and glue markings as follows:
	
	\begin{enumerate}
		\item[\textbf{(1)}] The markings $\tilde{t_{1}}$ and $ \tilde{t_{2}}$ on $\tilde{\mathbb{E}}$ are glued.
		
		\item[\textbf{(2)}] The marking $m_{i}$ on $\mathbb{E}$ is glued to the marking $\tilde{m_{i}}$ on $\tilde{\mathbb{E}}$, for each $i\in\mathbb{N}$.
		
		\item[\textbf{(3)}] The marking $m_{i}^{j}$ on $ \mathbb{E}$ is glued to the marking $\check{m}_{i}^{j}$ on $S(Id,h_j)$, for each $i\in\mathbb{N}$ and each $j\in \{1,\ldots, J\}$.
	\end{enumerate}
	Thus, we obtain the tame Loch Ness monster
	\begin{equation}\label{decorated_surface}
		S_{{\rm dec}},
	\end{equation}
	which is called \textbf{decorated surface} (see Figure \ref{Figure_decorated}).
	
	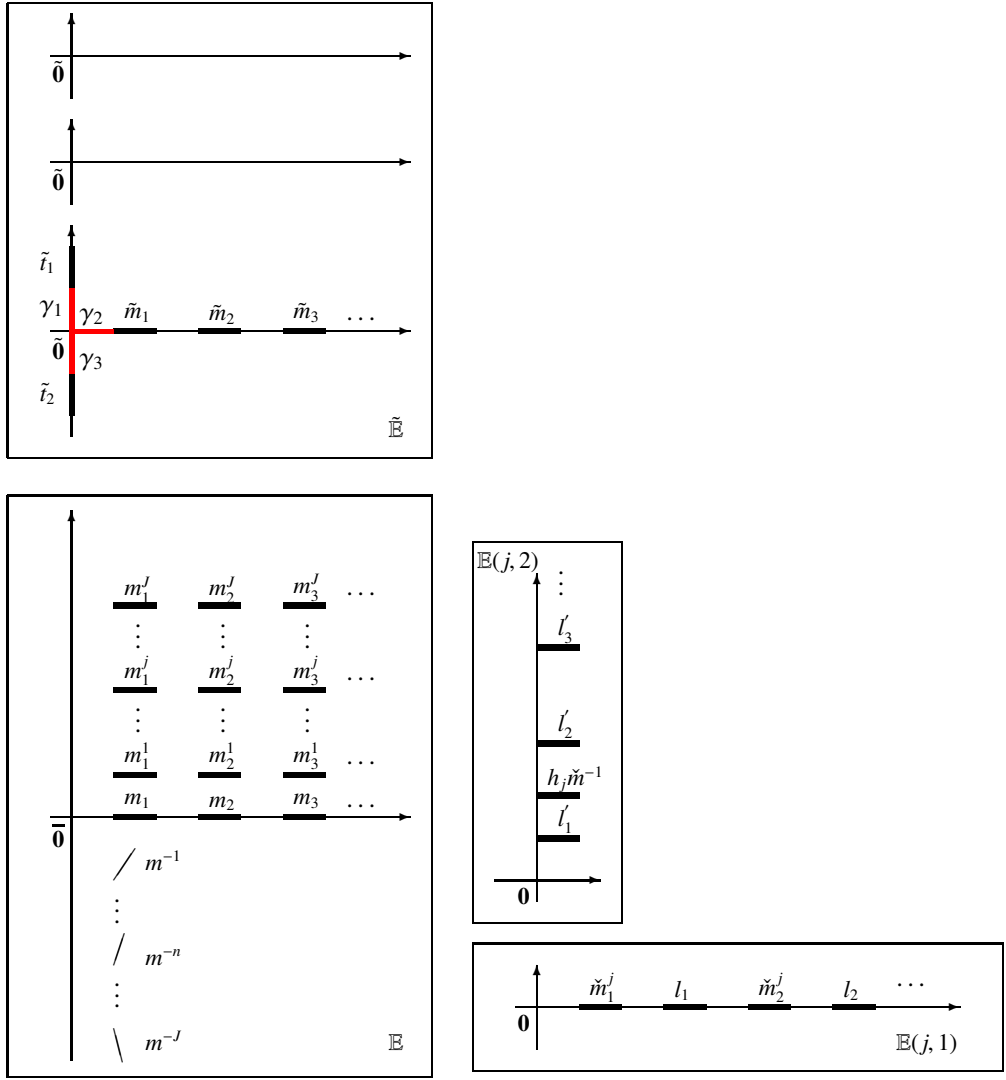
\begin{figure}[ht!]
		\begin{center}
			\setlength{\unitlength}{0.75pt}
			\begin{picture}(490,510)
				%\graphpaper(0,0)(480,530)
				\put(0,-3){\framebox(200,275)}
				\put(20,120){\vector(1,0){170}}
				\put(30,5){\vector(0,1){260}}
				\put(21,106){{\small$\overline{\textbf{0}}$}}
				%%%%%%%%%%%%%%marks_m
				\put(50,120){\linethickness{0.7mm}\line(1,0){20}}
				\put(55,126){{\small$m_{1}$}}
				\put(90,120){\linethickness{0.7mm}\line(1,0){20}}
				\put(95,125){\small{$m_{2}$}}
				\put(130,120){\linethickness{0.7mm}\line(1,0){20}}
				\put(135,126){{\small$m_{3}$}}
				\put(160,125){$\ldots$}
				\put(180,10){{\small$\mathbb{E}$}}
				%%%%%%%%%%%%%%%%marks_m^1
				\put(50,140){\linethickness{0.7mm}\line(1,0){20}}
				\put(55,146){{\small$m_{1}^{1}$}}
				\put(90,140){\linethickness{0.7mm}\line(1,0){20}}
				\put(95,146){{\small$m_{2}^{1}$}}
				\put(130,140){\linethickness{0.7mm}\line(1,0){20}}
				\put(135,146){{\small$m_{3}^{1}$}}
				\put(160,145){$\ldots$}
				%%%%%%%%%%%%%vdots
				\put(60,160){$\vdots$}
				\put(100,160){$\vdots$}
				\put(140,160){$\vdots$}
				%%%%%%%%%marks_m_j
				\put(50,180){\linethickness{0.7mm}\line(1,0){20}}
				\put(55,186){{\small$m_{1}^{j}$}}
				\put(90,180){\linethickness{0.7mm}\line(1,0){20}}
				\put(95,186){{\small$m_{2}^{j}$}}
				\put(130,180){\linethickness{0.7mm}\line(1,0){20}}
				\put(135,186){{\small$m_{3}^{j}$}}
				\put(160,185){$\ldots$}
				%%%%%%%%%%%%%%%%%%%%%%\vdots
				\put(60,200){$\vdots$}
				\put(100,200){$\vdots$}
				\put(140,200){$\vdots$}
				%%%%%%%%%%%%%%%%%%%%marks-m-J
				\put(50,220){\linethickness{0.7mm}\line(1,0){20}}
				\put(55,225){{\small$m_{1}^{J}$}}
				\put(90,220){\linethickness{0.7mm}\line(1,0){20}}
				\put(95,225){{\small$m_{2}^{J}$}}
				\put(130,220){\linethickness{0.7mm}\line(1,0){20}}
				\put(135,226){{\small $m_{3}^{J}$}}
				\put(160,225){$\ldots$}
				%%%%%%%%%%%%%%%%%%%%%% 
				\put(50,90){\linethickness{3mm}\line(2,3){10}} 
				\put(65,95){{\small$m^{-1}$}}
				\put(50,70){$\vdots$} 
				\put(50,50){\linethickness{3mm}\line(1,3){5}}
				\put(65,50){{\small $m^{-n}$}} 
				\put(50,30){$\vdots$}
				\put(50,20){\linethickness{3mm}\line(1,-4){4}}
				\put(65,10){{\small $m^{-J}$}} 
				%%%%%%%%%%%%%%%%%%%%%%%%%%%%%%%%%%%
				%%%%%%%%%%%%Covering
				\put(0,290){\framebox(200,215)}
				\put(20,350){\vector(1,0){170}}
				\put(30,300){\vector(0,1){100}}
				\put(30,310){\linethickness{0.7mm}\line(0,1){20}}
				\put(15,317){{\small$\tilde{t}_2$}}
				\put(30,370){\linethickness{0.7mm}\line(0,1){20}}
				\put(15,379){{\small$\tilde{t}_1$}}
				%%%%%%%%marks_on_covering
				\put(50,350){\linethickness{0.7mm}\line(1,0){20}}
				\put(30,350){{\color{red}\linethickness{0.5mm}\line(1,0){20}}}
				\put(30,350){{\color{red}\linethickness{0.6mm}\line(0,1){20}}}
				\put(30,350){{\color{red}\linethickness{0.6mm}\line(0,-1){20}}}
				\put(15,359){$\gamma_1$}
				\put(34,355){$\gamma_2$}
				\put(34,334){$\gamma_3$}
				\put(55,356){{\small$\tilde{m}_{1}$}}
				\put(90,350){\linethickness{0.7mm}\line(1,0){20}}
				\put(95,355){\small{$\tilde{m}_{2}$}}
				\put(130,350){\linethickness{0.7mm}\line(1,0){20}}
				\put(135,356){{\small$\tilde{m}_{3}$}}
				\put(160,355){$\ldots$}
				\put(180,300){{\small$\tilde{\mathbb{E}}$}}
				\put(21,337){{\small$\tilde{\textbf{0}}$}}
				\put(20,430){\vector(1,0){170}}
				\put(30,410){\vector(0,1){40}}
				\put(21,417){{\small$\tilde{\textbf{0}}$}}
				\put(20,480){\vector(1,0){170}}
				\put(30,460){\vector(0,1){40}}
				\put(21,468){{\small$\tilde{\textbf{0}}$}}
				
				%%%%%%%%%%%%%%%%%%%%%%%%%
				%%%%%%%%Buffer
				\put(220,70){\framebox(70,180)}
				\put(230,90){\vector(1,0){50}}
				\put(250,80){\vector(0,1){155}}
				\put(241,79){{\small$\textbf{0}$}}
				%%%%%%
				\put(250,110){\linethickness{0.8mm}\line(1,0){20}}
				\put(260,116){{\small$l_{1}^{'}$}}
				\put(250,130){\linethickness{0.8mm}\line(1,0){20}}
				\put(255,135){{\small$h_j\check{m}^{-1}$}}
				\put(250,155){\linethickness{0.8mm}\line(1,0){20}}
				\put(260,161){{\small$l_{2}^{'}$}}
				\put(250,200){\linethickness{0.8mm}\line(1,0){20}}
				\put(260,206){{\small$l_{3}^{'}$}}
				\put(260,225){$\vdots$}
				\put(222,238){{\small$\mathbb{E}(j,2)$}}
				%%%%%%%%%%%%%
				%%%%%%%%%%%%%%%%
				%%%%%%%%%
				\put(220,0){\framebox(232,60)}
				\put(240,30){\vector(1,0){210}}
				\put(250,10){\vector(0,1){40}}
				\put(241,19){{\small$\textbf{0}$}}
				%%%%%%%%%%%%%%
				\put(270,30){\linethickness{0.7mm}\line(1,0){20}}
				\put(275,36){{\small$\check{m}_{1}^{j}$}}
				\put(310,30){\linethickness{0.7mm}\line(1,0){20}}
				\put(315,35){{\small $l_{1}$}}
				\put(350,30){\linethickness{0.7mm}\line(1,0){20}}
				\put(355,36){{\small$\check{m}_{2}^{j}$}}
				\put(390,30){\linethickness{0.7mm}\line(1,0){20}}
				\put(395,35){{\small$l_{2}$}}
				\put(420,40){$\ldots$}
				\put(420,10){{\small$\mathbb{E}(j,1)$}}
			\end{picture}
			\caption{\emph{Decorated surface $S_{{\rm dec}}$.}}
			\label{Figure_decorated}
		\end{center}
	\end{figure} 
\end{construction}

\begin{remark}\label{r:singular_points_6pi}
	For each $j\in \{1,\ldots, J\}$, the markings  $h_{j}\check{m}^{-j}$ and $m^{-j}$ on the decorated surface $S_{{\rm dec}}$ have not yet been glued. Moreover, the surface $S_{{\rm dec}}$ has the following properties:
	
	\begin{enumerate} 
		\item[\textbf{(1)}] Its set of singular points consists of an infinitely many cone angle singularities of angle $4\pi$, and only one cone angle singularity of angle $6\pi$, which is denoted by $\tilde{\textbf{0}}$. 
		
		\item[\textbf{(2)}] There are only three saddle connections $\gamma_1$, $\gamma_2$, and $\gamma_3$, such that each one of them has as one of their endpoints the singularity $\tilde{\textbf{0}}$ (see Figure \ref{Figure_decorated}). The holonomy vectors of these saddle connections are $\{\pm e_1, \pm e_2\}$. 
	\end{enumerate}
	
	The surface $S_{{\rm dec}}$ is a slight modification of the surface appearing in Construction 4.6 in \cite{PSV}. In that construction, the authors introduced a tame Loch Ness monster with infinitely many markings on it. Nevertheless, in our case, we consider the same surface but with only a subset of these markings. Additionally, the decorated surfaces appearing in \cite{RamVal} cover different ends spaces; however, each of them has drawn an infinite family of markings for each element of $H$. This implies that our decorated surface $S_{{\rm dec}}$ is not studied in the aforementioned article. 	
\end{remark}

\subsubsection*{Step 2. The puzzle associated to the triplet $(1, G, H)$} Let $S_{g}$ be the affine copy of the decorated surface $S_{{\rm dec}}$, for each $g\in G$. We denote by $gh_{j}\check{m}^{-j}$ and $gm^{-j}$ (respectively) the markings on $S_g$, which are the images of the markings $h_{j}\check{m}^{-j}$ and $m^{-j}$ (respectively) via the affine diffeomorphism $f_{g}:S_{{\rm dec}}\to S_g$, where $j\in\{1,\ldots, J\}$. Thus, we define the \emph{puzzle associated to the triplet $(1,G,H)$} as 
\[
\mathfrak{P}(1,G,H):=\left\{S_{g}: g\in G\right\}.
\]
The following lemma will be used to prove the tameness of our surface $S$.

\begin{lemma}[\cite{PSV}]\label{lemma_PSV} 
	For every $g\in G$, the distance in $S_g$ between the families of markings $\{gh_{j}\check{m}^{-j}: j\in\{1,\ldots, J\} \}$ and $\{gm^{-j}: j\in\{1,\ldots, J\} \}$ is at least $1/\sqrt{2}$.
\end{lemma}

\subsubsection*{Step 3. The assembled surface $S$ to the puzzle $\mathfrak{P}(1,G,H)$} We consider the union $\bigcup_{g\in G} S_{g}$ and glue markings as follows: given the edge $(g,gh_{j})$ of the Cayley graph $Cay(G,H)$,
we glue the marking
$gh_{j}\check{m}^{-j}$ on $S_{g}$ to the marking $gh_{j}m^{-j}$ on $S_{gh_{j}}$.

We remark that, by construction, the markings $gh_{j}\check{m}^{-j}$ and  $gh_{j}m^{-j}$ are parallel, so the gluing is well-defined.  Thus, the \emph{assembled surface to the puzzle $\mathfrak{P}(1,G,H)$} obtained from the above gluing is a translation surface, which we denote by 
\begin{equation}
	S:=\bigcup\limits_{g\in G}S_{g}\bigg/ \sim.
\end{equation}

\subsection{The surface $S$ is a tame translation surface and its Veech group is the subgroup $G<{\rm GL}_{+}(2,\mathbb{R})$} 
One can use several of the ideas described in \cite[Theorem 3.7]{RamVal} to easily prove the following lemmas.

\begin{lemma}
	The translation surface $S$ is tame.
\end{lemma}	

\begin{proof}
	We must show that  $S$ is a complete metric space with respect to its natural flat metric $d$, and its set of singularities is discrete in $S$. Let $(\widehat{S},\widehat{d})$ be the metric completion space of $(S,d)$. For each $g\in G$,  we define the connected open subset
	\begin{equation}\label{eq:subsurface_removed}
		S'_{g}:=S_{g}\setminus \left\{ gh_{j}\check{m}^{-j}, gm^{-j}: j\in\{1,\ldots,J \}\right\}\subset S_{g},
	\end{equation}
	which is obtained from $S_{g}$ (see equation (\ref{decorated_surface})) by removing the markings  $gh_{j}\check{m}^{-j}$ and $gm^{-j}$ for each $j\in \{1,\ldots, J\}$. Using the inclusion map, the open subset $S'_{g}\subset S_{g}$ can be considered as a connected open subset of $S$. Then, the closure $\overline{S'_{g}}$ of $S'_{g}$ in $S$ is complete. If we take a Cauchy sequence $(x_n)_{n\in \mathbb{N}}$ in $S$ and the real number $\varepsilon=\frac{1}{2\sqrt{2}}$, then there is a positive integer $N(\varepsilon) \in \mathbb{N}$ such that for all natural numbers $m,n \geq N(\varepsilon)$, the terms $x_{m}, x_{n}$ satisfy  $\widehat{d}(x_{m},x_{n})< \varepsilon$. By Lemma \ref{lemma_PSV}, there is $g\in G$ such that the open ball $B_{\varepsilon}(x_{N(\varepsilon)})$ is contained in $\overline{S'_{g}}$. Since $\overline{B_{\varepsilon}(x_{N(\varepsilon)})}\subset \overline{S'(g)}$ is complete, the Cauchy sequence $(x_{n})_{n\in \mathbb{N}}$ converges within $\overline{B_{\varepsilon}(x_{N(\varepsilon)})}$. The discreteness of the singularities follows immediately from Lemma \ref{lemma_PSV}.
\end{proof}

\vspace{2mm}	
\begin{lemma}
	The Veech group of $S$ is $G$.
\end{lemma}

\begin{proof}
	Given that the group $G$ acts on $\mathfrak{P}(1,G,H):=\{S_g: g\in G\}$ by post-composition on charts,  then if we fix a matrix $\tilde{g}\in G$, for each $g\in G$, there exists a natural affine diffeomorphism $f_{\tilde{g}g}:S_{g} \to S_{\tilde{g}g}$, satisfying the following properties:
	
	\begin{enumerate}
		\item[\textbf{(1)}] The differential of $f_{\tilde{g}g}$ is the matrix $\tilde{g}$. 
		
		\item[\textbf{(2)}] The map $f_{\tilde{g}g}$ sends parallel markings to parallel markings. 
		
	\end{enumerate}
	Hence, the map $f: \bigcup\limits_{g\in G} S_{g} \to \bigcup\limits_{g\in G} S_{\tilde{g}g}$ defined by $f|_{S_{g}}:=f_{\tilde{g}g}$, is a gluing markings-preserving map. This yields an affine diffeomorphism in the quotient $F_{\tilde{g}}:S\to S$ with differential matrix $\tilde{g}$. Thus, we conclude that $G<\Gamma(S)$. Conversely, we consider $f:S\to S$ an affine orientation preserving diffeomorphism different from the identity. From Remark \ref{r:singular_points_6pi}, for each $g\in G$, the surface $S_{g}$ has one singularity of angle $6\pi$, which is denoted by $\tilde{\textbf{0}}_{g}$. There are only three saddle connections $\gamma_{1}^{g}$, $\gamma_{2}^{g}$, and $\gamma_{3}^{g}$ such that each one of them has that singularity as one of their endpoints. The holonomy vectors associated to these saddle connections are $\{\pm g\cdot e_1, \pm g\cdot e_2\}$. The function $f$ sends the singularity $\tilde{\textbf{0}}_{{\rm Id}}$ to the singularity $\tilde{\textbf{0}}_{g}$ for some $g\in G$, and  the differential matrix $df$ of $f$ must map $\{\pm e_1, \pm e_2\}$ to $\{\pm g\cdot e_1, \pm g\cdot e_2\}$. The only possibility is that $df=g$. Thus, we conclude that $\Gamma(S)<G$.
\end{proof}

\subsection{Ends space of the surface $S$} 

The description of the ends space of $S$, as stated in Theorem \ref{theorem.0.2}, follows from the following lemmas.

%%%%%%%%%%%%%%%%%%%%%%%%%%%%%%%%%%%%%%%%%%%%
%%%%%%%%%%%%%%%%%%%%%%%%%%%%%%%%%%%%%%%%%%% 	

\begin{lemma}\label{l:G_finite}
	If $G$ is finite, then the surface $S$ has as many ends as there are elements in the group $G$, and each end has infinite genus.
\end{lemma}

\begin{lemma}\label{l:G_no_finite}
	If $G$ is not finite, then the ends space of $S$  can be represented in the form  
	\[
	{\rm Ends}(S)={\rm Ends}_{\infty}(S)=\mathcal{B}\sqcup \mathcal{U},
	\]
	where $\mathcal{B}$ is a closed subset of ${\rm Ends}(S)$ homeomorphic to ${\rm Ends(G)}$, and $\mathcal{U}$ is a countable, dense, and open subset of ${\rm Ends}(S)$.
\end{lemma}

%%%%%%%%%%%%%%%%%%%%%%%%%%%%%%%%%%%%%%%%%
%%%%%%%%%%%%%%%%%%%%%%%%%%%%%%%%%%%%%%%%

\subsubsection*{Proof Lemma \ref{l:G_finite}} The group $G$ has cardinality $k$, for some $k\in\mathbb{N}$.  Let $K$ be a compact subset of $S$, we must prove that there exists a compact subset $K'\subset S$, such that $K\subset K'$, and $S\setminus K'$ consists of $k$ open connected components, and each one of them having infinite genus. 

For each $g\in G$, the affine copy $S_{g}$ is homeomorphic to the Loch Ness monster (see equation \eqref{decorated_surface}). Since the generating set $H$ of $G$ is finite, the set of markings 
\[
\left\{ gh_{j}\check{m}^{-j}, gm^{-j}: j\in\{1,\ldots, J\}\right\}
\]
on the affine copy $S_g$ is finite. We consider the connected subsurface $S'_{g}$ of $S_{g}$ as in equation \eqref{eq:subsurface_removed}, which has the following properties:

\begin{itemize} 
	\item[\textbf{(1)}] This subsurface $S'_{g}$ has infinite genus, and via the inclusion map, it can be considered as a connected subsurface of $S$ with infinite genus. 
	
	\item[\textbf{(2)}] The boundary $\partial S'_{g}$ of $S'_{g}$ in $S$ is compact because it is conformed by a finitely many disjoint closed curves. 
\end{itemize}
As $G$ is finite, from the preceding properties we hold that the set 
\[
S\setminus \bigcup\limits_{g\in G} \partial S'_{g}=\bigcup\limits_{g\in G} S'_{g},
\]
consists of $k$ open connected components, and each one of them having infinite genus.

On the other hand, let $K_{g}$ be the closure of the set $K\cap S'_{g}$ in $S_{g}$, for each $g\in G$. As $K_{g}$ is a compact subset of $S_{g}$, then there exists a compact subset $K'_{g}\subset S_{g}$ such that
\[
K_{g}\cup \left\{ gh_{j}\check{m}^{-j}, gm^{-j}: j\in\{1,\ldots, J\}\right\} \subset K'_{g},
\]
and $S_{g}\setminus K'_{g}$ consist of an open connected with infinite genus. We take $K'$ the closure of 
\[
\bigcup_{g\in G}\left(K'_{g}\setminus \{ gh_{j}\check{m}^{-j}, gm^{-j}: j\in\{1,\ldots, J\}\}\right)
\]
in $S$. As $G$ is finite, then $K'$ is a compact subset of $S$. By construction, we hold that $K\subset K'$, and the set
\[
S\setminus K'=\bigcup_{g\in G} (S_{g}\setminus K'_{g})\subset\bigcup\limits_{g\in G} S'_{g},
\]
consists of $k$ open connected components, and each one of them having infinite genus. \qed

%%%%%%%%%%%%%%%%%%%%%%%%%%%%%%%%%%%%%%%%%%%%%%%%%%%%%%%%
%%%%%%%%%%%%%%%%%%%%%%%%%%%%%%%
%%%%%%%%%%%%%%%%%%%%%G is not finite	

\subsubsection*{Proof Lemma \ref{l:G_no_finite}}  The sketch of the proof is the following. We begin by defining the set $\mathcal{U}$ from the ends of the affine copies $S_g$, and we will prove that it is a countable, discrete, and open subset of ${\rm Ends}(S)$. Then, we shall give an appropriate embedding $i_{\ast}$ from ${\rm Ends}(G)$ to ${\rm Ends}(S)$, where the image of ${\rm Ends}(G)$ under $i_{\ast}$ will be denoted by $\mathcal{B}$. By using an embedding from the Cayley graph ${\rm Cay}(G,H)$ to the surface $S$, we shall establish the equality 
\[
{\rm Ends}(S)={\rm Ends}_{\infty}(S)=\mathcal{B}\sqcup \mathcal{U},
\]
where $\mathcal{B}$ is closed, and $\mathcal{U}$ is an dense, and open subset of ${\rm Ends}(S)$.

%%%%%%%%%%%%%%%%%%%%%%%%%%%%%%%%%%%%%%%%%%%
%%%%%%%%%%%Step -1
%%%%%%%%%%%%%%%%%%%%%%%%%%%%%%%

\subsubsection*{Step 1. The set $\mathcal{U}$} For each $g\in G$, we take the subsurface $S'_{g}\subset S_{g}$ defined in equation \eqref{eq:subsurface_removed}. Recall that the boundary $\partial S'_{g}$ of the subsurface $S'_{g}$ is compact because it consists of a finitely many disjoint closed curves. Let $[U(g)_{n}]_{n\in\mathbb{N}}$ be the unique end of the Loch Ness monster $S_{g}$. Without loss of generality, we can assume that $U(g)_{n} \subset S'_{g}$ for each $n\in\mathbb{N}$. From the inclusion map, the surface $S'_{g}$ can be considered as a subsurface of $S$. Then the sequence $(U(g)_{n})_{n\in\mathbb{N}}$ of $S_{g}$ defines an end with infinite genus of the surface $S$.
\begin{remark}\label{r:properties_L}
	For any two different $g\neq \tilde{g}\in G$, the subsurfaces $S'_{g}$ and $S'_{\tilde{g}}$ of $S$ are disjoint.
\end{remark}

From the previous Remark, we obtain the countable set $\mathcal{U}$ conformed by different ends of $S$, given by
\begin{equation}\label{set_U}
	\mathcal{U}:=\left\{[U(g)_{n}]_{n\in\mathbb{N}}\in {\rm Ends}(S):g\in G\right\}\subset {\rm Ends}(S).
\end{equation}

Let us note that the subset $\mathcal{U}\subset {\rm Ends}(S)$ is both discrete and open. This is a consequence of the following fact. For each  $g\in G$, the open subset $U(g)_{1}$ of $S$ has a compact boundary $\partial U(g)_{1}$ in $S$. Thus, we define the open subset $(U(g)_{1})^{\ast}$ of ${\rm Ends}(S)$, which satisfies 
\[
(U(g)_{1})^{\ast}\cap \mathcal{U}=\left\{[U(g)_{n}]_{n\in\mathbb{N}}\right\}.
\] 

%%%%%%%%%%%%%%%%%%%%%%%%%%%%%%%%
%%%%%%%%%%Step -2
%%%%%%%%%%%%%%%%%%%%%%

\subsubsection*{Step 2. The embedding $i_{\ast}:{\rm Ends}(G)\hookrightarrow {\rm Ends}(S)$} Let $\overline{S'_{g}}$ be the closure in $S$ of the surface $S'_{g}$ (see equation \eqref{eq:subsurface_removed}). Given a non-empty connected open subset $W$ of ${\rm Cay}(G,H)$ with compact boundary $\partial W$, we can, suppose without loss of generality, that the boundary $\partial W \subset V({\rm Cay}(G,H))=G$. We then define the subset $\tilde{W}\subset S$ given by 
\begin{equation}\label{eq:induced_open}
	\tilde{W}:={\rm Int} \left( \bigcup_{g\in G \cap (W\cup\partial W)} \overline{S'_{g}} \right)\subset S.
\end{equation}
This set $\tilde{W}$ is a non-empty, connected, and open subset of $S$ with a compact boundary. Moreover, it is a subsurface of $S$ with infinite genus. In the following Remark, we state two properties of this object, which can be easily deduced.

\begin{remark}\label{remark_intersection}
	Given $W$ and  $V$ two non-empty, connected, and open subsets of ${\rm Cay}(G,H)$ each one them having compact boundary $\partial W$ and $\partial V$, respectively, such that $\partial W, \partial V \subset G$, then
	\begin{itemize}	
		\item[\textbf{(1)}] If $W\supset V$, then $\tilde{W}\supset\tilde{V}$.
		
		\item[\textbf{(2)}] If $W\cap V= \emptyset$, then $\tilde{W} \cap\tilde{V}=\emptyset$.
	\end{itemize}
\end{remark}

From the above Remark, the end $\left[W_{n}\right]_{n\in \mathbb{N}}$ of the group $G$ naturally defines the end $[\tilde{W}_{n}]_{n\in\mathbb{N}}$ of the surface $S$, which has infinite genus. Hence, we obtain a well-define map $i_{\ast}:{\rm Ends}(G)\to  {\rm Ends}(S)$, given by
\begin{equation}\label{eq:embedding_EndG}
	[W_{n}]_{n\in \mathbb{N}} \mapsto [\tilde{W}_{n}]_{n\in \mathbb{N}}.
\end{equation}

\begin{claim}
	The map $i_{\ast}$ is an embedding.
\end{claim}

\begin{proof}
	We must show that $i_{\ast}$ is \emph{injective}. Let $[W_{n}]_{n\in\mathbb{N}}$ and $[V_{n}]_{n\in\mathbb{N}}$ be two different ends of $G$. Then, there is $l\in \mathbb{N}$ such that $W_{l}\cap V_{l}=\emptyset$. By item (2) of Remark \ref{remark_intersection}, it follows that $\tilde{W}_{l}\cap \tilde{V}_{l}=\emptyset$. It proves that the ends $i_{\ast}([W_{n}]_{n\in\mathbb{N}})=[\tilde{W}_{n}]_{n\in\mathbb{N}}$ and  $i_{\ast}([V_{n}]_{n\in\mathbb{N}})=[\tilde{V}_{n}]_{n\in\mathbb{N}}$ in ${\rm Ends}(S)$ are different.
	
	\emph{Continuity}. We consider an end $[W_{n}]_{n\in\mathbb{N}}$ of the group $G$ and an open subset $V\subset S$ with a compact boundary, such that $i_{\ast}([W_{n}]_{n\in\mathbb{N}})=[\tilde{W}_{n}]_{n\in\mathbb{N}} \in V^{\ast}\subset {\rm Ends}(S)$. We must prove that there is a neighborhood $Z^{\ast}\subset {\rm Ends}(G)$ of $[W_{n}]_{n\in\mathbb{N}}$ such that $i_{\ast}\left(Z^{\ast}\right)\subset V^{\ast}$. Given that $[\tilde{W}_{n}]_{n\in\mathbb{N}} \in V^{\ast}$, there exists some $k\in\mathbb{N}$ such that
	\begin{equation}\label{eq:contained_W_j}
		\tilde{W}_{k} \subset V,
	\end{equation}
	We take the open subset $W_{k}$ of the Cayley graph ${\rm Cay}(G,H)$, which defines the open subset $\tilde{W}$ (see equation \eqref{eq:induced_open}), and consider the open
	\[
	Z^{\ast}:=(W_{k})^{\ast},
	\]
	of ${\rm Ends}(G)$, which is a neighborhood  of $[W_{n}]_{n\in\mathbb{N}}$. To ensure that $i_{\ast}(Z^{\ast})\subset V^{\ast}$, we consider any end $[U_{n}]_{n\in\mathbb{N}}\in {\rm Ends}(G)$ such that $[U_{n}]_{n\in\mathbb{N}}\in Z^{\ast}=(W_{k})^{\ast}$, and check that $i_{\ast}([U_{n}]_{n\in\mathbb{N}})=[\tilde{U}_{n}]_{n\in\mathbb{N}}\in V^{\ast}$. Since $U_{m}\subset W_{k}$ for some $m\in\mathbb{N}$, it follows from item (1) of Remark \ref{remark_intersection} that
	\[
	\tilde{U}_{m}\subset \tilde{W}_{k}.
	\]
	As $\tilde{W}_{k}\subset V$, we conclude that $\tilde{U}_{m}\subset V$, which implies that $i_{\ast}([U_{n}]_{n\in\mathbb{N}})=[\tilde{U}_{n}]_{n\in\mathbb{N}}\in V^{\ast}$.
	
	Finally, the map $i_{\ast}$ is \emph{closed} because any continuous map from a compact space to a Hausdorff space is closed. Therefore, $i_{\ast}$ is an embedding.
\end{proof}

We denote the image of the map $i_{\ast}$ as 
\[
\mathcal{B}:=i_{\ast}({\rm Ends}(G)).
\]
From the definition of the set $\mathcal{U}$ given in equation \eqref{set_U}, we conclude that $\mathcal{B}\cap \mathcal{U}=\emptyset$, and $\mathcal{B}\sqcup \mathcal{U}\subset {\rm Ends}(S)$.

%%%%%%%%%%%%%%%%%%%%%%%%%%%%%%%%%%%%%%%%%%%%%%%%%%%%%%%%%%%%%
%%%%%%%%%%%%%%%%%%%%%%embedding i
%%%%%%%%%%%%%%%%%%%%%%%%%

\subsubsection*{Step 3. The embedding $i:{\rm Cay}(G,H)\hookrightarrow S$} We now describe the image of each vertex and edge of ${\rm Cay}(G,H)$ under the map $i$. 

For each $g\in G$, let $\overline{\textbf{0}}_{g}$ denote the point in the affine copy $S_{g}$ that corresponds to the image of the point $\overline{\textbf{0}}$ (see equation \eqref{decorated_surface}) in the decorated surface $S_{{\rm  dec}}$ via the affine diffeomorphism $f_{g}:S_{{\rm dec}}\to S_{g}$. Then the surface $S'_{g}$ described in equation \eqref{eq:subsurface_removed}, contains the point $\overline{\textbf{0}}_{g}$. Thus, we define the map $h:  V({\rm Cay}(G,H))=G \to  S$ given by  
\begin{equation}\label{eq:map_vertices} %(ec:3.43)
	g  \mapsto \overline{\textbf{0}}_{g}.
\end{equation}

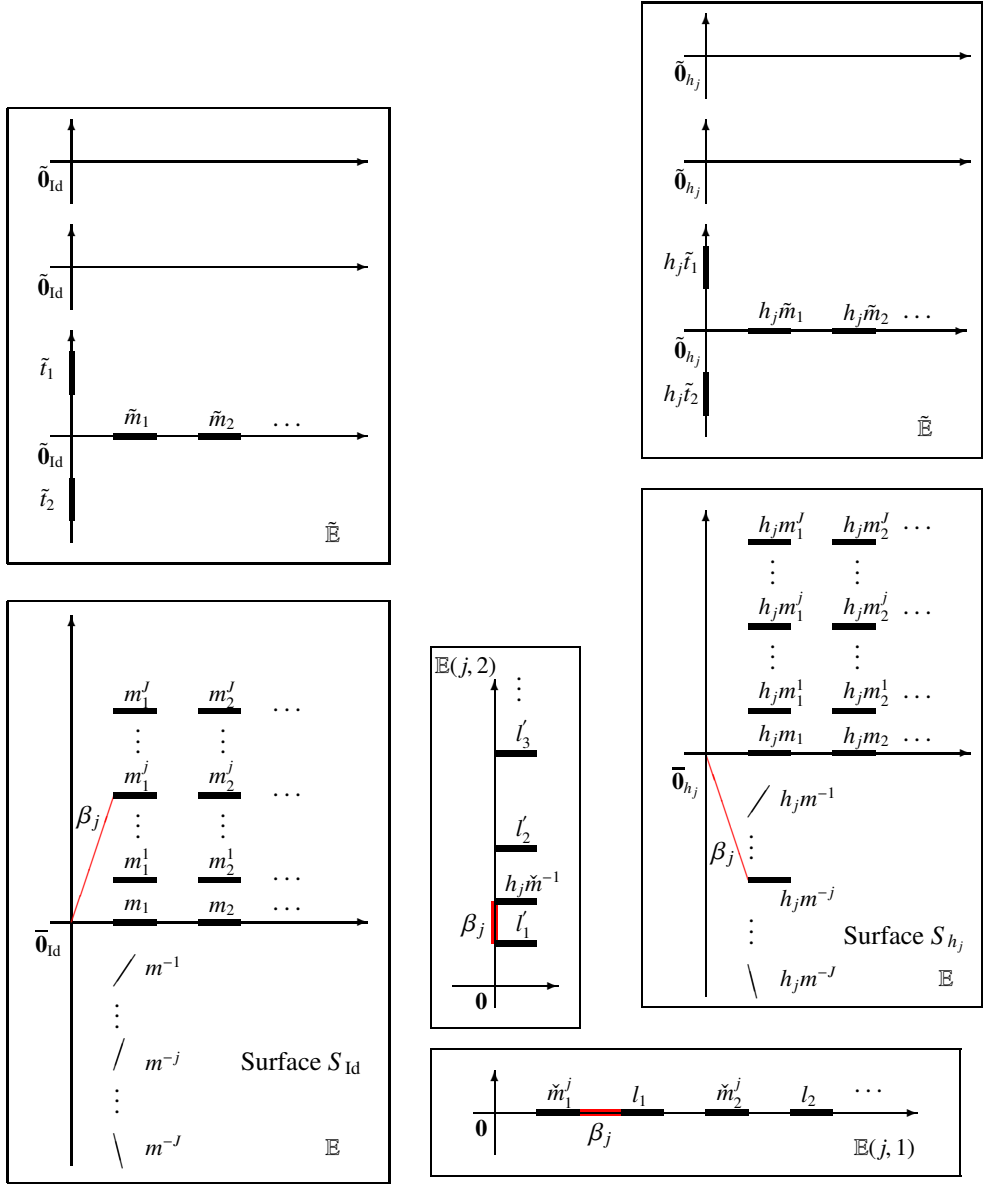
\begin{figure}[ht!]
	\begin{center}
		\setlength{\unitlength}{0.8pt}
		\begin{picture}(490,540)
			%\graphpaper(0,0)(480,530)
			\put(0,-3){\framebox(180,275)}
			\put(20,120){\vector(1,0){150}}
			\put(30,5){\vector(0,1){260}}
			\put(13,106){{\small$\overline{\textbf{0}}_{\rm Id}$}}
			\put(30,120){{\color{red}\linethickness{1.8mm}\line(1,3){20}}}
			\put(33,167){$\beta_j$}
			\put(215,116){$\beta_j$}
			\put(275,17){$\beta_j$}
			\put(333,150){$\beta_j$}
			\put(230,110){{\color{red}\linethickness{0.7mm}\line(0,1){20}}}
			\put(270,30){{\color{red}\linethickness{0.7mm}\line(1,0){20}}}
			\put(350,140){{\color{red}\linethickness{1.8mm}\line(-1,3){20}}}
			%%%%%%%%%%%%%%marks_m
			\put(50,120){\linethickness{0.7mm}\line(1,0){20}}
			\put(55,126){{\small$m_{1}$}}
			\put(90,120){\linethickness{0.7mm}\line(1,0){20}}
			\put(95,125){\small{$m_{2}$}}
			\put(125,126){$\ldots$}
			\put(150,10){{\small$\mathbb{E}$}}
			\put(110,50){Surface $S_{\rm Id}$}
			%%%%%%%%%%%%%%%%marks_m^1
			\put(50,140){\linethickness{0.7mm}\line(1,0){20}}
			\put(55,146){{\small$m_{1}^{1}$}}
			\put(90,140){\linethickness{0.7mm}\line(1,0){20}}
			\put(95,146){{\small$m_{2}^{1}$}}
			\put(125,140){$\ldots$}
			%%%%%%%%%%%%%vdots
			\put(60,160){$\vdots$}
			\put(100,160){$\vdots$}
			\put(50,180){\linethickness{0.7mm}\line(1,0){20}}
			\put(55,186){{\small$m_{1}^{j}$}}
			\put(90,180){\linethickness{0.7mm}\line(1,0){20}}
			\put(95,186){{\small$m_{2}^{j}$}}
			\put(125,180){$\ldots$}
			%%%%%%%%%%%%%%%%%%%%%%\vdots
			\put(60,200){$\vdots$}
			\put(100,200){$\vdots$}
			%%%%%%%%%%%%%%%%%%%%marks-m-J
			\put(50,220){\linethickness{0.7mm}\line(1,0){20}}
			\put(55,225){{\small$m_{1}^{J}$}}
			\put(90,220){\linethickness{0.7mm}\line(1,0){20}}
			\put(95,225){{\small$m_{2}^{J}$}}
			\put(125,220){$\ldots$}
			%%%%%%%%%%%%%%%%%%%%%% 
			\put(50,90){\linethickness{0.7mm}\line(2,3){10}} 
			\put(65,95){{\small$m^{-1}$}}
			\put(50,70){$\vdots$} 
			\put(50,50){\linethickness{3mm}\line(1,3){5}}
			\put(65,50){{\small $m^{-j}$}} 
			\put(50,30){$\vdots$}
			\put(50,20){\linethickness{3mm}\line(1,-4){4}}
			\put(65,10){{\small $m^{-J}$}} 
			%%%%%%%%%%%%%%%%%%%%%%%%%%%%%%%%%%%
			%%%%%%%%%%%%Covering
			\put(0,290){\framebox(180,215)}
			\put(20,350){\vector(1,0){150}}
			\put(30,300){\vector(0,1){100}}
			\put(30,310){\linethickness{0.7mm}\line(0,1){20}}
			\put(15,317){{\small$\tilde{t}_2$}}
			\put(30,370){\linethickness{0.7mm}\line(0,1){20}}
			\put(15,379){{\small$\tilde{t}_1$}}
			%%%%%%%%marks_on_covering
			\put(50,350){\linethickness{0.7mm}\line(1,0){20}}
			\put(55,356){{\small$\tilde{m}_{1}$}}
			\put(90,350){\linethickness{0.7mm}\line(1,0){20}}
			\put(95,355){\small{$\tilde{m}_{2}$}}
			\put(125,355){$\ldots$}
			\put(150,300){{\small$\tilde{\mathbb{E}}$}}
			\put(14,337){{\small$\tilde{\textbf{0}}_{{\rm Id}}$}}
			\put(20,430){\vector(1,0){150}}
			\put(30,410){\vector(0,1){40}}
			\put(14,417){{\small$\tilde{\textbf{0}}_{{\rm Id}}$}}
			\put(20,480){\vector(1,0){150}}
			\put(30,460){\vector(0,1){40}}
			\put(14,468){{\small$\tilde{\textbf{0}}_{\rm Id}$}}
			
			%%%%%%%%%%%%%%%%%%%%%%%%%
			%%%%%%%%Buffer
			\put(200,70){\framebox(70,180)}
			\put(210,90){\vector(1,0){50}}
			\put(230,80){\vector(0,1){155}}
			\put(221,79){{\small$\textbf{0}$}}
			%%%%%%
			\put(230,110){\linethickness{0.8mm}\line(1,0){20}}
			\put(240,116){{\small$l_{1}^{'}$}}
			\put(230,130){\linethickness{0.8mm}\line(1,0){20}}
			\put(235,135){{\small$h_j\check{m}^{-1}$}}
			\put(230,155){\linethickness{0.8mm}\line(1,0){20}}
			\put(240,161){{\small$l_{2}^{'}$}}
			\put(230,200){\linethickness{0.8mm}\line(1,0){20}}
			\put(240,206){{\small$l_{3}^{'}$}}
			\put(240,225){$\vdots$}
			\put(202,238){{\small$\mathbb{E}(j,2)$}}
			%%%%%%%%%%%%%
			%%%%%%%%%%%%%%%%
			%%%%%%%%%
			\put(200,0){\framebox(250,60)}
			\put(220,30){\vector(1,0){210}}
			\put(230,10){\vector(0,1){40}}
			\put(221,19){{\small$\textbf{0}$}}
			%%%%%%%%%%%%%%
			\put(250,30){\linethickness{0.7mm}\line(1,0){20}}
			\put(255,36){{\small$\check{m}_{1}^{j}$}}
			\put(290,30){\linethickness{0.7mm}\line(1,0){20}}
			\put(295,35){{\small $l_{1}$}}
			\put(330,30){\linethickness{0.7mm}\line(1,0){20}}
			\put(335,36){{\small$\check{m}_{2}^{j}$}}
			\put(370,30){\linethickness{0.7mm}\line(1,0){20}}
			\put(375,35){{\small$l_{2}$}}
			\put(400,40){$\ldots$}
			\put(400,10){{\small$\mathbb{E}(j,1)$}}
			%%%%%%%%%%%%%%%%%%%%%%%%%%%
			%%%%%%%%%%%%%%%%%%%%%%%%%%%%%%%
			%%%%%%%Decorated- S_h
			%%%%%%%%%%%%%%%%%%%%
			%%%%%%%%%%%%%%%%%%
			\put(300,80){\framebox(160,245)}
			\put(320,200){\vector(1,0){135}}
			\put(330,85){\vector(0,1){230}}
			\put(314,183){{\small$\overline{\textbf{0}}_{h_{j}}$}}
			%%%%%%%%%%%%%%marks_m
			\put(350,200){\linethickness{0.7mm}\line(1,0){20}}
			\put(355,206){{\small$h_j m_{1}$}}
			\put(390,200){\linethickness{0.7mm}\line(1,0){20}}
			\put(395,205){\small{$h_j m_{2}$}}
			\put(423,206){$\ldots$}
			\put(395,110){Surface $S_{h_j}$}
			\put(440,90){{\small$\mathbb{E}$}}
			%%%%%%%%%%%%%%%%marks_m^1
			\put(350,220){\linethickness{0.7mm}\line(1,0){20}}
			\put(355,226){{\small$h_j m_{1}^{1}$}}
			\put(390,220){\linethickness{0.7mm}\line(1,0){20}}
			\put(395,226){{\small$h_j m_{2}^{1}$}}
			\put(423,226){$\ldots$}
			%%%%%%%%%%%%%vdots
			\put(360,240){$\vdots$}
			\put(400,240){$\vdots$}
			%%%%%%%%%marks_m_j
			\put(350,260){\linethickness{0.7mm}\line(1,0){20}}
			\put(355,266){{\small$h_j m_{1}^{j}$}}
			\put(390,260){\linethickness{0.7mm}\line(1,0){20}}
			\put(395,266){{\small$h_j m_{2}^{j}$}}
			\put(423,266){$\ldots$}
			%%%%%%%%%%%%%%%%%%%%%%\vdots
			\put(360,280){$\vdots$}
			\put(400,280){$\vdots$}
			%%%%%%%%%%%%%%%%%%%%marks-m-J
			\put(350,300){\linethickness{0.7mm}\line(1,0){20}}
			\put(355,305){{\small$h_j m_{1}^{J}$}}
			\put(390,300){\linethickness{0.7mm}\line(1,0){20}}
			\put(395,305){{\small$h_j m_{2}^{J}$}}
			\put(423,306){$\ldots$}
			%%%%%%%%%%%%%%%%%%%%%% 
			\put(350,170){\linethickness{3mm}\line(2,3){10}} 
			\put(365,175){{\small$h_j m^{-1}$}}
			\put(350,150){$\vdots$} 
			\put(350,140){\linethickness{0.7mm}\line(1,0){20}}
			\put(365,128){{\small $h_j m^{-j}$}} 
			\put(350,110){$\vdots$}
			\put(350,100){\linethickness{3mm}\line(1,-4){4}}
			\put(365,90){{\small $h_j m^{-J}$}} 
			%%%%%%%%%%%%%%%%%%%%%%%%%%%%%%%%%%%
			%%%%%%%%%%%%Covering
			\put(300,340){\framebox(160,215)}
			\put(320,400){\vector(1,0){133}}
			\put(330,350){\vector(0,1){100}}
			\put(330,360){\linethickness{0.7mm}\line(0,1){20}}
			\put(310,367){{\small$h_j \tilde{t}_2$}}
			\put(330,420){\linethickness{0.7mm}\line(0,1){20}}
			\put(310,429){{\small$h_j \tilde{t}_1$}}
			%%%%%%%%marks_on_covering
			\put(350,400){\linethickness{0.7mm}\line(1,0){20}}
			\put(355,406){{\small$h_j \tilde{m}_{1}$}}
			\put(390,400){\linethickness{0.7mm}\line(1,0){20}}
			\put(395,405){\small{$h_j \tilde{m}_{2}$}}
			\put(423,406){$\ldots$}
			\put(430,350){{\small$\tilde{\mathbb{E}}$}}
			\put(315,387){{\small$\tilde{\textbf{0}}_{h_j}$}}
			\put(320,480){\vector(1,0){135}}
			\put(330,460){\vector(0,1){40}}
			\put(315,467){{\small$\tilde{\textbf{0}}_{h_j}$}}
			\put(320,530){\vector(1,0){135}}
			\put(330,510){\vector(0,1){40}}
			\put(315,518){{\small$\tilde{\textbf{0}}_{h_j}$}}
		\end{picture}
		\caption{\emph{Image of $\beta_j$.}}
		\label{Figure_beta}
	\end{center}
\end{figure}

On the other hand, for each $j\in\{1,\ldots, J\}$, there is a simple polygonal path $\beta_{j}:[0,1]\to S$ satisfying the following properties:
\begin{itemize}
	\item[\textbf{(1)}] The initial and terminal points of $\beta_{j}$ are  $\overline{\textbf{0}}_{\rm Id}$ and $\overline{\textbf{0}}_{h_{j}}$, respectively. See Figure \ref{Figure_beta}.
	
	\item[\textbf{(2)}] For each $i\neq j \in\{1,\ldots, J\}$,  the intersection $\beta_{i}([0,1])\cap\beta_{j}([0,1])=\left\{\overline{\textbf{0}}_{\rm Id}\right\}$. 
\end{itemize}
%%%%%%%%%%%%%%%%%%%%%%%%%%%%%
%%%%%%%%%%%%%%%
%%%%%%%%%%%%%%%%%%
%%%%%%%%%%%%
Since the edge $({\rm Id},h_j)$ of the Cayley graph ${\rm Cay}(G,H)$ is homeomorphic to the open interval $(0,1)$, we can, suppose without loss of generality, assume that the curve $\beta_{j}$ is defined from $[Id,h_j]$ to $S$ such that  $\beta_j({\rm Id})=\overline{\textbf{0}}_{\rm Id}$ and $\beta_j(h_j)=\overline{\textbf{0}}_{h_j}$. Given that the Veech group of the surface $S$ is $G$, for each $g\in G$, there is an affine diffeomorphism $f_g: S\to S$ whose differential is $d f_{g}=g$. Thus, we get the composition path
\begin{equation}\label{eq:composition_curve}
	f_g \circ\beta_j:[0,1]\to S,
\end{equation}
satisfying the following properties:
\begin{itemize}
	\item[\textbf{(1)}] The initial and terminal points of $f_g\circ \beta_{j}$ are  $\overline{\textbf{0}}_{g}$ and $\overline{\textbf{0}}_{g h_{j}}$, respectively.
	
	\item[\textbf{(2)}] For each $i\neq j \in\{1,\ldots, J\}$,  the intersection $f_g\circ \beta_{i}([0,1])\cap f_g\circ \beta_{j}([0,1])=\left\{\overline{\textbf{0}}_{g}\right\}$. 
\end{itemize}
Similarly, since the edge $(g,gh_j)$ of the Cayley graph ${\rm Cay}(G,H)$ is homeomorphic to the open interval $(0,1)$, we can, suppose without loss of generality, assume that the composition path $f_g \circ\beta_j$ is defined from $[g,gh_{j}]$ to $S$ such that $f_g \circ \beta_j(g)=\overline{\textbf{0}}_{g}$ and $f_g\circ \beta_j(gh_j)=\overline{\textbf{0}}_{gh_j}$.

From equations \eqref{eq:map_vertices} and \eqref{eq:composition_curve}, we obtain the embedding 
\begin{equation}\label{eq:embedding_i}
	i : {\rm Cay}(G,H) \hookrightarrow  S,
\end{equation} 
such that $i_{|G}:=h$ and $i_{|[g,gh_j]}:=f_{g} \circ \beta_{j}$ for each $g\in G$ and $j\in\{1,\ldots,J\}$.

%%%%%%%%%%%%%%%%%%%%%%%%%%%%%%%%%%%%%%%%%%%%
%%%%%%%%%%%%%%Step-4
%%%%%%%%%%%%%%%%%%%%%

\subsubsection*{Step 4. The equality ${\rm Ends}(S)=\mathcal{B}\sqcup \mathcal{U}$} We must prove that $ {\rm Ends}(S)\subset  \mathcal{B}\sqcup \mathcal{U}$.  Let $[U_n]_{n\in\mathbb{N}}$ be an end of $S$. Since $S=\bigcup\limits_{g\in G}\overline{S'_{g}}$, for each $n\in \mathbb{N}$, we consider the subset 
\[
G(n)=\left\{g\in G: \overline{S'_{g}}\cap U_{n}\neq \emptyset\right\}\subset G,
\]
and we define the open subset
\begin{equation*}
	Z_{n}:={\rm Int}\left(\bigcup\limits_{g\in G(n)} \overline{S'_{g}}\right)\subset S,
\end{equation*}
which has the following properties:

\begin{itemize}
	\item[\textbf{(1)}] Since $U_{n}$ is an non-empty, connected, and open subset of $S$ with compact boundary, the set $Z_{n}$ is also a connected, and open subset of $S$ with compact boundary, for each $n\in\mathbb{N}$.
	
	\item[\textbf{(2)}] As $U_{n}\supset U_{n+1}$, it follows that $Z_{n}\supset Z_{n+1}$ for each $n\in\mathbb{N}$. 
\end{itemize}

Using the definition of end and the construction of $Z_n$, it is easy to show that the sequences $(Z_{n})_{n\in \mathbb{N}}$ and $(U_{n})_{n\in\mathbb{N}}$ defines the same end of $S$. In other words, $[U_{n}]_{n\in\mathbb{N}}=[Z_{n}]_{n\in\mathbb{N}}$. We shall now prove that the end $[Z_{n}]_{n\in\mathbb{N}}$ belongs to $\mathcal{B}\sqcup\mathcal{U}$. We notice that one of the following cases must occur:

\textbf{Case 1.} There is $N\in\mathbb{N}$ such that $G(N)$ is finite. then there exists $g\in G$ such that for all $m\geq N$, we hold  
\[
Z_{m} \subset S'_{g}.
\]
This implies that the sequences $(Z_{n})_{n\in \mathbb{N}}$ and $(U(g)_{n})_{n\in\mathbb{N}}$ must be equivalent (see equation \eqref{set_U}). Thus, $[U_n]_{n\in\mathbb{N}}\in \mathcal{U}$.

\textbf{Case 2.} Otherwise, for each $n\in\mathbb{N}$, the subset $G(n)\subset G$ is infinite. As the embedding $i$, described in equation \eqref{eq:embedding_i}, is a continuous map, the inverse image  
\[
{\hat Z}_{n}:=i^{-1}\left(Z_{n}\cap i({\rm Cay}(G,H))\right),
\]
is a connected and open subset of ${\rm Cay}(G,H)$ with compact boundary for each $n\in\mathbb{N}$. Moreover, the sequence $(\hat{Z}_{n})_{n\in\mathbb{N}}$ defines an end of the group $G$. By the construction of the sequence $(Z_{n})_{n\in\mathbb{N}}$ of $S$, the embedding $i_{\ast}$, defined in  \eqref{eq:embedding_EndG}, sends the end $[\hat{Z}_{n}]_{n\in \mathbb{N}}$ of $G$ to the end $[Z_{n}]_{n\in \mathbb{N}}$ of $S$. This implies that $[Z_{n}]_{n\in\mathbb{N}}$ belongs to $\mathcal{B}$. Thus, we conclude that ${\rm Ends}(S)=\mathcal{B}\sqcup \mathcal{U}$.

\subsubsection*{Step 5. The set $\mathcal{B}$ is closed and the set $\mathcal{U}$ is a dense open} Since $\mathcal{U}$ is an open subset of ${\rm Ends}(S)$, its complement ${\rm Ends}(S)\setminus \mathcal{U}=\mathcal{B}$ is a closed subset of ${\rm Ends}(S)$. We shall prove that $\mathcal{U}$ is dense. Let $[Z_{n}]_{n\in \mathbb{N}}$ be an end of $\mathcal{B}$, we must show that this end belongs to the closure of $\mathcal{U}$. 

Let $U$ be a non-empty, connected, and open subset of $S$ with compact boundary such that the open subset $U^{\ast}\subset {\rm Ends }(S)$ contains the end $[Z_{n}]_{n\in\mathbb{N}}$. There exists $\tilde{g}\in \{ g\in G: \overline{S'_{g}}\cap U\neq \emptyset\}$ such that $S'_{\tilde{g}}\subset U$. This condition implies that the end $[U(\tilde{g})_{n}]_{n\in\mathbb{N}}$ of $\mathcal{U}$ belongs to $U^{\ast}$. Therefore, the end $[Z_{n}]_{n\in\mathbb{N}}$ is in the closure of $\mathcal{U}$.
\qed

%%%%%%%%%%%%%%%%%%%%%%%%%%%%%%%%%%%%%
%%%%%%%%%%%%%%%%%%%%%%%%%%%%%%%%%%%

%%%%%%%%%%%%%%%%%%%%%%%%%%%%%%%%%%%%%%%%%%%
%%%%%%%%%%%%%%%%%%%%%%%%%%%%%%%%%%%%%%%%%%%%%%%

\subsection*{Acknowledgements}
Camilo Ram\'irez Maluendas expresses his gratitude to the UNIVERSIDAD NACIONAL DE COLOMBIA, SEDE MANIZALES.  He has dedicated this work to his beautiful family: Marbella and Emilio, in appreciation of their love and support.

%%%%%%%%%%%%%%%%%%%%%%%%%%%%%%%%%%%%%%%%%%%%%%%%%%%%%%%%%%%%%%%%%%%%%%%%%%%%%%%%%%%%%%%%%%%%%%%%%%%%%%%%%%%%%%%%%%%%%%%%%%%%%%%%%%%%%%%%%%%%%%%%%%%%%%%%%%%%%%%%%%

\end{document}